\documentclass[11pt]{article}
\usepackage{amsmath,amssymb,amstext,dsfont,fancyvrb,float,fontenc,graphicx,subfigure, theorem}

\title{\Large\bf On Hilbert and Riemann problems.\\ An alternative approach.}

\author{\sc Vladimir Ryazanov}
\date{}

\cleardoublepage 
\usepackage[strict]{changepage}
\def\abstractname{Abstract -}   
\def\abstract{\begin{adjustwidth}{1cm}{1cm} \par    \footnotesize \noindent {\bf \abstractname}
\def\endabstract{ \end{adjustwidth} \smallskip }}

{\theorembodyfont{\itshape}}
{\theorembodyfont{\itshape}}
{\theorembodyfont{\itshape}}
{\theorembodyfont{\itshape}}
{\theorembodyfont{\itshape}}
{\theorembodyfont{\rm}}
{\theorembodyfont{\rm}}
{\theorembodyfont{\rm}}
{\theorembodyfont{\rm }}

\begin{document}
\maketitle
\vskip 1.5em

\vskip 1.5em
 \begin{abstract}
Recall that the Hilbert (Riemann-Hilbert) boundary value problem was
recently solved  in \cite{R1} for arbitrary measurable coefficients
and for arbitrary measurable boundary data in terms of nontangential
limits and principal asymptotic values. Here it is developed a new
approach making possible to obtain new results on tangential limits.
It is shown that the spaces of the found solutions have the infinite
dimension for prescribed collections of Jordan arcs terminating in
almost every boundary point. Similar results are proved for the
Riemann problem.
 \end{abstract}


\begin{MSC}
primary   31A05, 31A20, 31A25, 31B25, 35Q15; se\-con\-da\-ry 30E25,
31C05, 34M50, 35F45.
\end{MSC}

\bigskip

\section{Introduction}

The Hilbert (Riemann-Hilbert) boundary value problem, the Riemann
and Poincare boundary value problems are basic in the theory of
analytic functions and they are closely interconnected, see for the
history e.g. the monographs \cite{G}, \cite{M} and \cite{V}, and
also the last works \cite{R1}-\cite{R4}.

\medskip

Recall that the classical setting of the {\bf Riemann problem} in a
smooth Jordan domain $D$ of the complex plane $\mathbb{C}$ was on
finding analytic functions $f^+: D\to\mathbb C$ and $f^-:\mathbb
C\setminus \overline{D}\to\mathbb C$ that admit continuous
extensions to $\partial D$ and satisfy the boundary condition
\begin{equation}\label{eqRIEMANN} f^+(\zeta)\ =\
A(\zeta)\cdot f^-(\zeta)\ +\ B(\zeta) \quad\quad\quad \forall\
\zeta\in\partial D
\end{equation}
with prescribed H\"older continuous functions $A: \partial
D\to\mathbb C$ and $B: \partial D\to\mathbb C$.

\medskip

Recall also that the {\bf Riemann problem with shift} in $D$ was on
finding such functions $f^+: D\to\mathbb C$ and $f^-:\mathbb
C\setminus \overline{D}\to\mathbb C$ satisfying the condition
\begin{equation}\label{eqSHIFT} f^+(\alpha(\zeta))\ =\ A(\zeta)\cdot
f^-(\zeta)\ +\ B(\zeta) \quad\quad\quad \forall\ \zeta\in\partial D
\end{equation}
where $\alpha :\partial D\to\partial D$ was a one-to-one sense
preserving correspondence having the non-vanishing H\"older
continuous derivative with respect to the natural parameter on
$\partial D$. The function $\alpha$ is called a {\bf shift
function}. The special case $A\equiv 1$ gives the so--called {\bf
jump problem}.

\medskip

The classical setting of the {\bf Hilbert (Riemann-Hilbert) boundary
value problem} was on finding analytic functions $f$ in a domain
$D\subset\mathbb C$ bounded by a rectifiable Jordan curve with the
boundary condition
\begin{equation}\label{eqHILBERT}
\lim\limits_{z\to\zeta}\ \mathrm {Re}\
\{\overline{\lambda(\zeta)}\cdot f(z)\}\ =\ \varphi(\zeta)
\quad\quad\quad\ \ \ \forall \ \zeta\in \partial D
\end{equation}
with functions $\lambda$ and $\varphi$ that are continuously
differentiable with respect to the natural parameter $s$ on
$\partial D$ and, moreover, $|\lambda|\ne 0$ everywhere on $\partial
D$. Hence without loss of generality one can assume that
$|\lambda|\equiv 1$ on $\partial D$.

\medskip

It is clear that if we start to consider the Hilbert and Riemann
problems with measurable boundary data, the requests on the
existence of the limits at all points $\zeta\in\partial D$ and along
all paths terminating in $\zeta$ lose any sense (as we{ll as the
conception of the index). Thus, the notion of solutions of the
Hilbert and Riemann problems should be widened in this case. The
nontangential limits were a suitable tool from the function theory
of one complex variable, see e.g. \cite{R1}--\cite{R4}. Here it is
proposed an alternative approach admitting tangential limits.

\bigskip

Given a Jordan curve $C$ in $\mathbb C$, we say that a family of
Jordan arcs $\{J_{\zeta}\}_{\zeta\in C}$ is of {\bf class \cal{BS}
(of Bagemihl--Seidel class)}, cf. \cite{BS}, 740--741, if all
$J_{\zeta}$ lie in a ring $\frak{R}$ generated by $C$ and a Jordan
curve $C_*$ in $\mathbb C$, $C_*\cap C={\O}$, $J_{\zeta}$ is joining
$C_*$ and $\zeta\in C$, every $z\in \frak{R}$ belongs to a single
arc $J_{\zeta}$, and for a sequence of mutually disjoint Jordan
curves $C_n$ in $\frak{R}$ such that $C_n\to C$ as $n\to\infty$,
$J_{\zeta}\cap C_n$ consists of a single point for each $\zeta\in C$
and $n=1,2, \ldots$.

\bigskip

In particular, a family of Jordan arcs $\{J_{\zeta}\}_{\zeta\in C}$
is of class \cal{BS} if $J_{\zeta}$ is generated by an isotopy of
$C$. For instance, every curvilinear ring $\frak{R}$ one of whose
boundary component is $C$ can be mapped with a conformal mapping $g$
onto a circular ring $R$ and the inverse mapping
$g^{-1}:R\to\frak{R}$ maps radial lines in $R$ onto suitable Jordan
arcs $J_{\zeta}$ and centered circles in $R$ onto Jordan curves
giving the corresponding isotopy of $C$ to other boundary component
of $\frak{R}$.

\bigskip

Now, if $\Omega\subset\mathbb C$ is an open set bounded by a finite
collection of mutually disjoint Jordan curves, then we say that a
family of Jordan arcs $\{J_{\zeta}\}_{\zeta\in \partial\Omega}$ is
of class \cal{BS} if its restriction to each component of
$\partial\Omega$ is so.

\bigskip

\section{On the Hilbert problem}

{\bf Theorem 1.} {\it\, Let $D$ be a bounded domain in $\mathbb C$
whose boundary consists of a finite number of mutually disjoint
rectifiable Jordan curves, and let $\lambda:\partial D\to\mathbb C$,
$|\lambda (\zeta)|\equiv 1$, $\varphi:\partial D\to\mathbb R$ and
$\psi:\partial D\to \mathbb R$ be measurable functions with respect
to the natural parameter. Suppose that $\{
\gamma_{\zeta}\}_{\zeta\in\partial D}$ is a family of Jordan arcs of
class \cal{BS} in ${D}$.

\medskip

Then there exist single-valued analytic functions $f: D\to\mathbb C$
such that
\begin{equation}\label{eqRE} \lim\limits_{z\to\zeta}\ \mathrm {Re}\
\{\overline{\lambda(\zeta)}\cdot f(z)\}\ =\ \varphi(\zeta)\ ,
\end{equation}
\begin{equation}\label{eqIM}
\lim\limits_{z\to\zeta}\ \mathrm {Im}\
\{\overline{\lambda(\zeta)}\cdot f(z)\}\ =\ \psi(\zeta)\
\end{equation}
along $\gamma_{\zeta}$ for a.e. $\zeta\in\partial D$ with respect to
the natural parameter.}

\medskip

{\bf Remak 1.} Thus, the space of all solutions $f$ of the Hilbert
problem (\ref{eqRE}) in the given sense has the infinite dimension
for any prescribed $\varphi$, $\lambda$ and $\{
\gamma_{\zeta}\}_{\zeta\in D}$ because the space of all measurable
functions $\psi:\partial D\to \mathbb R$ has the infinite dimension.

\medskip

\begin{proof} Indeed, set $\Psi(\zeta)=\varphi(\zeta)+i\psi(\zeta)$
and $\Phi(\zeta)=\lambda(\zeta)\cdot\Psi(\zeta)$ for all
$\zeta\in\partial D$. Then by Theorem 2 in \cite{BS} there is a
single-valued analytic function $f$ such that
\begin{equation}\label{eqBS} \lim\limits_{z\to\zeta}\ f(z)\ =\ \Phi(\zeta)
\end{equation}
along $\gamma_{\zeta}$ for a.e. $\zeta\in\partial D$ with respect to
the natural parameter. Then also \begin{equation}\label{eqLIM}
\lim\limits_{z\to\zeta}\ \overline{\lambda(\zeta)}\cdot f(z)\ =\
\Psi(\zeta)
\end{equation}
along $\gamma_{\zeta}$ for a.e. $\zeta\in\partial D$ with respect to
the natural parameter.
\end{proof}

\medskip

Similar result can be formulated for arbitrary Jordan domains in
terms of the harmonic measure.

\medskip

{\bf Theorem 2.} {\it\, Let $D$ be a bounded domain in $\mathbb C$
whose boundary consists of a finite number of mutually disjoint
Jordan curves, and let $\lambda:\partial D\to\mathbb C$, $|\lambda
(\zeta)|\equiv 1$, $\varphi:\partial D\to\mathbb R$ and
$\psi:\partial D\to \mathbb R$ be measurable functions with respect
to the harmonic measure. Suppose that $\{
\gamma_{\zeta}\}_{\zeta\in\partial D}$ is a family of Jordan arcs of
class \cal{BS} in ${D}$.

\medskip

Then there exist single-valued analytic functions $f: D\to\mathbb C$
such that
\begin{equation}\label{eqRE-A} \lim\limits_{z\to\zeta}\ \mathrm {Re}\
\{\overline{\lambda(\zeta)}\cdot f(z)\}\ =\ \varphi(\zeta)\ ,
\end{equation}
\begin{equation}\label{eqIM-A}
\lim\limits_{z\to\zeta}\ \mathrm {Im}\
\{\overline{\lambda(\zeta)}\cdot f(z)\}\ =\ \psi(\zeta)\
\end{equation}
along $\gamma_{\zeta}$ for a.e. $\zeta\in\partial D$ with respect to
the harmonic measure.}

\bigskip

{\bf Remak 2.} Again, the space of all solutions $f$ of the
Riemann-Hilbert problem (\ref{eqRE-A}) in the given sense has the
infinite dimension for any prescribed $\varphi$, $\lambda$ and $\{
\gamma_{\zeta}\}_{\zeta\in D}$ because the space of all functions
$\psi:\partial D\to \mathbb R$ that are measurable with respect to
the harmonic measure has the infinite dimension.

\bigskip

\begin{proof} Theorem 2 is reduced to Theorem 1 in the following way.

\medskip

First, there is a conformal mapping $\omega$ of $D$ onto a circular
domain $\mathbb D_*$  whose boundary consists of a finite number of
circles and points, see e.g. Theorem V.6.2 in \cite{Go}. Note that
$\mathbb D_*$ is not degenerate because isolated singularities of
conformal mappings are removable that is due to the well-known
Weierstrass theorem, see e.g. Theorem 1.2 in \cite{CL}. Applying in
the case of need the inversion with respect to a boundary circle of
$\mathbb D_*$, we may assume that $\mathbb D_*$ is bounded.

\medskip

Remark that $\omega$ is extended to a homeomorphism $\omega_*$ of
$\overline D$ onto $\overline{\mathbb D_*}$, see e.g. point (i) of
Lemma 3.1 in \cite{R4}. Set $\Lambda = \lambda\circ\Omega$, $\Phi =
\varphi\circ\Omega$ and $\Psi = \psi\circ\Omega$ where $\Omega :
\partial{\mathbb D_*}\to\partial D$ is the restriction of $\Omega_* :=\omega_*^{-1}$ to $\partial{\mathbb D_*}$.
Let us show that these functions are measurable with respect to the
natural parameter on $\partial{\mathbb D_*}$.

\medskip

For this goal, note first of all that the sets of the harmonic
measure zero are invariant under conformal mappings between multiply
connected Jordan domains because a composition of a harmonic
function with a conformal mapping is again a harmonic function.
Moreover, a set $E\subset
\partial{\mathbb D_*}$ has the harmonic measure zero if and only if
it has the length zero, say in view of the integral representation
of the harmonic measure through the Green function of the domain
$\mathbb D_*$, see e.g. Section II.4 in \cite{N}.

\medskip

Hence $\Omega$ and $\Omega^{-1}$ transform measurable sets into
measurable sets because every measurable set is the union of a
sigma-compact set and a set of measure zero, see e.g. Theorem
III(6.6) in \cite{S}, and continuous mappings transform compact sets
into compact sets. Thus, the functions $\lambda$, $\varphi$ and
$\psi$ are measurable with respect to the harmonic measure on
$\partial D$ if and only if the functions $\Lambda$, $\Phi$ and
$\Psi$ are measurable with respect to the natural parameter on
$\partial\mathbb D_*$.

\medskip

Then by Theorem 1 there exist single-valued analytic functions $F:
D\to\mathbb C$ such that
\begin{equation}\label{eqRE-R} \lim\limits_{w\to\xi}\ \mathrm {Re}\
\{\overline{\Lambda(\xi)}\cdot F(w)\}\ =\ \Phi(\xi)\ ,
\end{equation}
\begin{equation}\label{eqIM-R}
\lim\limits_{w\to\xi}\ \mathrm {Im}\ \{\overline{\Lambda(\xi)}\cdot
F(w)\}\ =\ \Psi(\xi)\
\end{equation}
along $\Gamma_{\xi}=\omega(\gamma_{\Omega(\xi)})$ for a.e.
$\xi\in\partial \mathbb D_*$ with respect to the natural parameter.

\medskip

Thus, by the construction the functions $f=F\circ\omega$ are the
desired analytic functions $f:D\to\mathbb C$ satisfying the boundary
conditions (\ref{eqRE-A}) and (\ref{eqIM-A}) along $\gamma_{\zeta}$
for a.e. $\zeta\in\partial D$ with respect to the harmonic measure.
\end{proof}

\bigskip

{\bf Remark 3.} Many investigations were devoted to the nonlinear
Hilbert (Riemann-Hilbert) boundary value problems with conditions of
the type
\begin{equation}\label{eqNONLINEAR} \Phi(\,\zeta,\, f(\zeta)\, )\ =\ 0 \quad\quad\quad \forall\
\zeta\in\partial D\ ,
\end{equation}
see e.g. \cite{EW}, \cite{K} and \cite{W}. It is natural also to
weaken such conditions to
\begin{equation}\label{eqNONLINEAR} \Phi(\,\zeta,\, f(\zeta)\, )\ =\ 0 \quad\quad\quad \mbox{for a.e.}\quad
\zeta\in\partial D\ .
\end{equation}
It is easy to see that the proposed approach makes possible  also to
reduce such problems to the algebraic and measurable solvability of
the relation
\begin{equation}\label{eqNONLINEAR} \Phi(\,\zeta,\, v\, )\ =\ 0
\end{equation}
with respect to a complex-valued function $v(\zeta)$, cf. e.g.
\cite{Gr}.


Through suitable modifications of $\Phi$ under the corresponding
mappings of Jordan boundary curves onto the unit circle $\mathbb S =
\{ \zeta\in\mathbb C : |\zeta|=1\}$, we may assume that $\zeta$
belongs to $\mathbb S\mathbf{}$.


\section{On the Riemann problem}

{\bf Theorem 3.} {\it\, Let $D$ be a domain in $\overline{\mathbb
C}$ whose boundary consists of a finite number of mutually disjoint
rectifiable Jordan curves, $A: \partial D\to\mathbb C$ and $B:
\partial D\to\mathbb C$ be measurable functions with respect to the
natural parameter. Suppose that
$\{\gamma^+_{\zeta}\}_{\zeta\in\partial D}$ and
$\{\gamma^-_{\zeta}\}_{\zeta\in\partial D}$ are families of Jordan
arcs of class \cal{BS} in ${D}$ and $\mathbb C\setminus\overline{
D}$, correspondingly.


Then there exist single-valued analytic functions $f^+: D\to\mathbb
C$ and $f^-:\overline{\mathbb C}\setminus\overline{D}\to\mathbb C$
that satisfy (\ref{eqRIEMANN}) for a.e. $\zeta\in\partial D$ with
respect to the natural parameter where $f^+(\zeta)$ and $f^-(\zeta)$
are limits of $f^+(z)$ and $f^-(z)$ az $z\to\zeta$ along
$\gamma^+_{\zeta}$ and $ \gamma^-_{\zeta}$, correspondingly.


Furthermore, the space of all such couples $(f^+,f^-)$ has the
infinite dimension for every couple $(A, B)$ and any collections
$\gamma^+_{\zeta}$ and $ \gamma^-_{\zeta}$, $\zeta\in\partial D$.}

\medskip

Theorem 3 is a special case of the following lemma on the
generalized Riemann problem with shifts that can be useful for other
goals, too.

\medskip

{\bf Lemma 1.} {\it\, Under the hypotheses of Theorem 3, let in
addition $\alpha : \partial D\to\partial D$ be a homeomorphism
keeping components of $\partial D$ such that $\alpha$ and
$\alpha^{-1}$ have the $(N)-$property of Lusin with respect to the
natural parameter.


Then there exist single-valued analytic functions $f^+: D\to\mathbb
C$ and $f^-:\overline{\mathbb C}\setminus\overline{D}\to\mathbb C$
that satisfy (\ref{eqSHIFT}) for a.e. $\zeta\in\partial D$ with
respect to the natural parameter where $f^+(\zeta)$ and $f^-(\zeta)$
are limits of $f^+(z)$ and $f^-(z)$ az $z\to\zeta$ along
$\gamma^+_{\zeta}$ and $ \gamma^-_{\zeta}$, correspondingly.


Furthermore, the space of all such couples $(f^+,f^-)$ has the
infinite dimension for every couple $(A, B)$ and any collections
$\gamma^+_{\zeta}$ and $ \gamma^-_{\zeta}$, $\zeta\in\partial D$.}

\bigskip

\begin{proof} First, let $D$ be bounded and let $g^-: \partial D\to\mathbb C$
be a measurable function. Note that the function
\begin{equation}\label{eqCONNECTION}g^+\ :=\ \{ A\cdot
g^- + B\}\circ \alpha^{-1}\end{equation} is measurable. Indeed,
$E:=\{ A\cdot g^- + B\}^{-1}(\Omega)$ is a measurable subset of
$\partial D$ for every open set $\Omega\subseteq\mathbb C$ because
the function $A\cdot g^- + B$ is measurable by the hypotheses. Hence
the set $E$ is the union of a sigma-compact set and a set of measure
zero, see e.g. Theorem III(6.6) in \cite{S}. However, continuous
mappings transform compact sets into compact sets and, thus,
$\alpha(E)=\alpha\circ\{ A\cdot g^- + B\}
^{-1}(\Omega)=(g^+)^{-1}(\Omega)$ is a measurable set, i.e. the
function $g^+$ is really measurable.

Then by Theorem 2 in \cite{BS} there is a single-valued analytic
function $f^+: D\to\mathbb C$ such that
\begin{equation}\label{eqBS} \lim\limits_{z\to\xi}\ f^+(z)\ =\ g^+(\xi)
\end{equation}
along $\gamma^+_{\xi}$ for a.e. $\xi\in\partial D$ with respect to
the natural parameter. Note that $g^+(\alpha(\zeta))$ is determined
by the given limit for a.e. $\zeta\in\partial D$ because
$\alpha^{-1}$ also has the $(N)-$property of Lusin.

Note that $\overline{\mathbb C}\setminus\overline{D}$ consists of a
finite number of (simply connected) Jordan domains $D_0, D_1, \ldots
, D_m$ in the extended complex plane $\overline{\mathbb C}=\mathbb
C\cup \{ \infty\}$. Let $\infty\in D_0$. Then again by Theorem 2 in
\cite{BS} there exist single-valued analytic functions $f_l^-:
D_l\to\mathbb C$, $l=1, \ldots , m,$ such that
\begin{equation}\label{eqBS} \lim\limits_{z\to\zeta}\ f_l^-(z)\ =\
g_l^-(\zeta)\ ,\quad\quad\quad g_l^-:=g^-|_{\partial D_l}\ ,
\end{equation}
along $\gamma^-_{\zeta}$ for a.e. $\zeta\in\partial D_l$ with
respect to the natural parameter.

Now, let $S$ be a circle that contains $D$ and let $j$ be the
inversion of $\overline{\mathbb C}$ with respect to $S$. Set
$$D_*=j(D_0),\quad g_*=\overline {g_0\circ j},\quad  g_0^-:=g^-|_{\partial
D_0},\quad \gamma^*_{\xi}=j\left(\gamma^-_{j(\xi)}\right),\quad
\xi\in\partial D_*\ .$$ Then by Theorem 2 in \cite{BS} there is a
single-valued analytic function $f_*: D_*\to\mathbb C$ such that
\begin{equation}\label{eqBS} \lim\limits_{w\to\xi}\ f_*(w)\ =\ g_*(\xi)
\end{equation}
along $\gamma^*_{\xi}$ for a.e. $\xi\in\partial D_*$ with respect to
the natural parameter. Note that $f_0^-:=\overline {g_*\circ j}$ is
a single-valued analytic function in $D_0$ and by construction
\begin{equation}\label{eqBS0} \lim\limits_{z\to\zeta}\ f_0^-(z)\ =\
g_0^-(\zeta)\ ,\quad\quad\quad g_0^-:=g^-|_{\partial D_0}\ ,
\end{equation} along $\gamma^-_{\zeta}$ for a.e. $\zeta\in\partial D_0$ with
respect to the natural parameter.

Thus, the functions $f^-_l$, $l=0, 1, \ldots , m,$ form an analytic
function $f^-: \overline{\mathbb C}\setminus\overline{D}\to\mathbb
C$ satisfying (\ref{eqSHIFT}) for a.e. $\zeta\in\partial D$ with
respect to the natural parameter.

The space of all such couples $(f^+,f^-)$ has the infinite dimension
for every couple $(A, B)$ and any collections $\gamma^+_{\zeta}$ and
$ \gamma^-_{\zeta}$, $\zeta\in\partial D$, in view of the above
construction because of the space of all measurable functions $g^-:
\partial D\to\mathbb C$ has the infinite dimension.

The case of unbounded $D$ is reduced to the case of bounded $D$
through the complex conjugation and the inversion of
$\overline{\mathbb C}$ with respect to a circle $S$ in some of the
components of $\overline{\mathbb C}\setminus\overline{D}$ arguing as
above.
\end{proof}

\bigskip

{\bf Remark 4.} Some investigations were devoted also to the
nonlinear Riemann problems with boundary conditions of the form
\begin{equation}\label{eqNONLINEAR} \Phi(\,\zeta,\, f^+(\zeta),\, f^-(\zeta)\, )\ =\ 0 \quad\quad\quad \forall\
\zeta\in\partial D\ .
\end{equation}
It is natural as above to weaken such conditions to the following
\begin{equation}\label{eqNONLINEAR} \Phi(\,\zeta,\, f^+(\zeta),\, f^-(\zeta)\, )\ =\ 0 \quad\quad\quad \mbox{for a.e.}\quad
\zeta\in\partial D\ .
\end{equation}
It is easy to see that the proposed approach makes possible  also to
reduce such problems to the algebraic and measurable solvability of
the relations
\begin{equation}\label{eqNONLINEAR} \Phi(\,\zeta,\, v,\, w\, )\ =\ 0
\end{equation}
with respect to complex-valued functions $v(\zeta)$ and $w(\zeta)$,
cf. e.g. \cite{Gr}.

\medskip

Through suitable modifications of $\Phi$ under the corresponding
mappings of Jordan boundary curves onto the unit circle $\mathbb S =
\{ \zeta\in\mathbb C : |\zeta|=1\}$, we may assume that $\zeta$
belongs to $\mathbb S\mathbf{}$.

\bigskip

{\bf Example 1.} For instance, correspondingly to the scheme given
above, special nonlinear problems of the form
\begin{equation}\label{eqNONLINEAR} f^+(\zeta)\ =\ \varphi(\,\zeta,\,
 f^-(\zeta)\, ) \quad\quad\quad \mbox{for
a.e.}\quad \zeta\in\partial\mathbb D
\end{equation}
in the unit disk $\mathbb D=\{ z\in\mathbb C: |z|<1\}$ are always
solved if the function $\varphi : \mathbb S\times\mathbb C\to\mathbb
C$ satisfies the {\bf Caratheodory conditions}: $\varphi(\zeta, w)$
is continuous in the variable $w\in\mathbb C$ for a.e.
$\zeta\in\mathbb S$ and it is measurable in the variable
$\zeta\in\mathbb S$ for all $w\in\mathbb C$.

\medskip

Furthermore, the spaces of solutions of such problems always have
the infinite dimension. Indeed, the function
$\varphi(\zeta,\psi(\zeta))$ is measurable in $\zeta\in\mathbb S$
for every measurable function $\psi:\mathbb S\to\mathbb C$ if the
function $\varphi$ satisfies the {Caratheodory conditions}, see e.g.
Section 17.1 in \cite{KZPS}, and the space of all measurable
functions $\psi:\mathbb S\to\mathbb C$ has the infinite dimension.

\bigskip

{\bf Problems.} Finally, it is necessary to point out the open
problems on solvability of Hilbert and Riemann problems along any
prescribed families of arcs but not only along families of the
Bagemihl--Seidel class and, more generally, along any prescribed
families of paths to a.e. boundary point.




 { \footnotesize
\medskip
\medskip
 \vspace*{1mm}

\noindent {\it Vladimir Ryazanov}\\
Institute of Applied Mathematics and Mechanics \\
National Academy of Sciences of Ukraine \\
1 Dobrovol'skogo Str., 84100 Slavyansk, UKRAINE\\
E-mail: {\tt vl.ryazanov1@gmail.com}}

 \end{document}